\pdfoutput=1
\documentclass{amsart}
\usepackage{geometry}
\usepackage{graphicx}
\usepackage[colorlinks,citecolor=blue,urlcolor=blue]{hyperref}

\usepackage{amssymb, amsmath, amsthm}
\newcommand{\eq}{\begin{equation}}
\newcommand{\en}{\end{equation}}

\newtheorem{Theorem}{Theorem}
\newtheorem{theorem}[Theorem]{Theorem}
\newtheorem{lemma}[Theorem]{Lemma}
\newtheorem{corollary}[Theorem]{Corollary}
\newtheorem{construction}[Theorem]{Construction}
\newtheorem{proposition}[Theorem]{Proposition}
\newtheorem{example}[Theorem]{Example}
\newtheorem{defn}[Theorem]{Definition}

\newtheorem{question}[Theorem]{Question}
\newtheorem{conjecture}[Theorem]{Conjecture}
\newtheorem{condition}[Theorem]{Condition}
\newtheorem{problem}[Theorem]{Problem}
\theoremstyle{definition}
\newtheorem{remark}[Theorem]{Remark}

\def\cmin {convex minorant}

\newcommand {\unC} {\underline{C}}
\newcommand {\unc} {\underline{c}}

\newcommand {\IR} {\mbox{\msbm\symbol{'122}}}
\newcommand{\IP}{\mathbb{P}}
\newcommand{\IE}{\mathbb{E}}

\newcommand{\te}{\rightarrow}
\newcommand{\ed}{\mbox{$ \ \stackrel{d}{=}$ }}

\newcommand {\ee}{e}

\def\IR{\mathbb{R}}

\def\ba#1\ee{\begin{align*}#1\end{align*}}
\def\ban#1\ee{\begin{align}#1\end{align}}

\newcommand {\argmin} {{\rm argmin}}

\newcommand{\E}{\IE}

\newcommand{\tr}{$(\tau, \rho)$}
\newtheorem{definition}[Theorem]{Definition}

\newcommand{\son}{S^{[0,n]}}

\begin{document}

\title{Convex minorants of random walks and L\'evy processes}

\author{Josh Abramson}
\address{University of California, Berkeley}
\email{\emph{JA}: josh@stat.Berkeley.EDU, \emph{JP}: pitman@stat.Berkeley.EDU, \emph{NR}: ross@stat.Berkeley.EDU}
\author{Jim Pitman}
\thanks{JP's research supported in part by N.S.F.\ Grant DMS-0806118} 
\author{Nathan Ross}
\author{Ger\'onimo Uribe Bravo}
\address{Universidad Nacional Aut\'onoma de M\'exico}\thanks{GUB's research supported by a postdoctoral fellowship from UC MexUS - CoNaCyt and N.S.F. Grant DMS-0806118}\email{geronimo@matem.unam.mx}
\begin{abstract} 
This article provides an overview of recent work on descriptions and properties of the \cmin\ of random walks and L\'evy processes
as detailed in \cite{ap10,  piro10, pub10}, which summarize and extend the literature on these subjects.

The results surveyed include point process descriptions of the convex minorant of random walks and L\'evy processes on a fixed finite interval, up to an independent exponential time, and in the infinite horizon case.  These descriptions follow from the invariance of these processes under an adequate path transformation.  In the case of Brownian motion, we note how further special properties of this process, including time-inversion, imply a sequential description for the convex minorant of the Brownian meander. 
\end{abstract}

\maketitle

\section{Introduction}

\begin{figure}
  \centering
  \includegraphics[width=\textwidth]{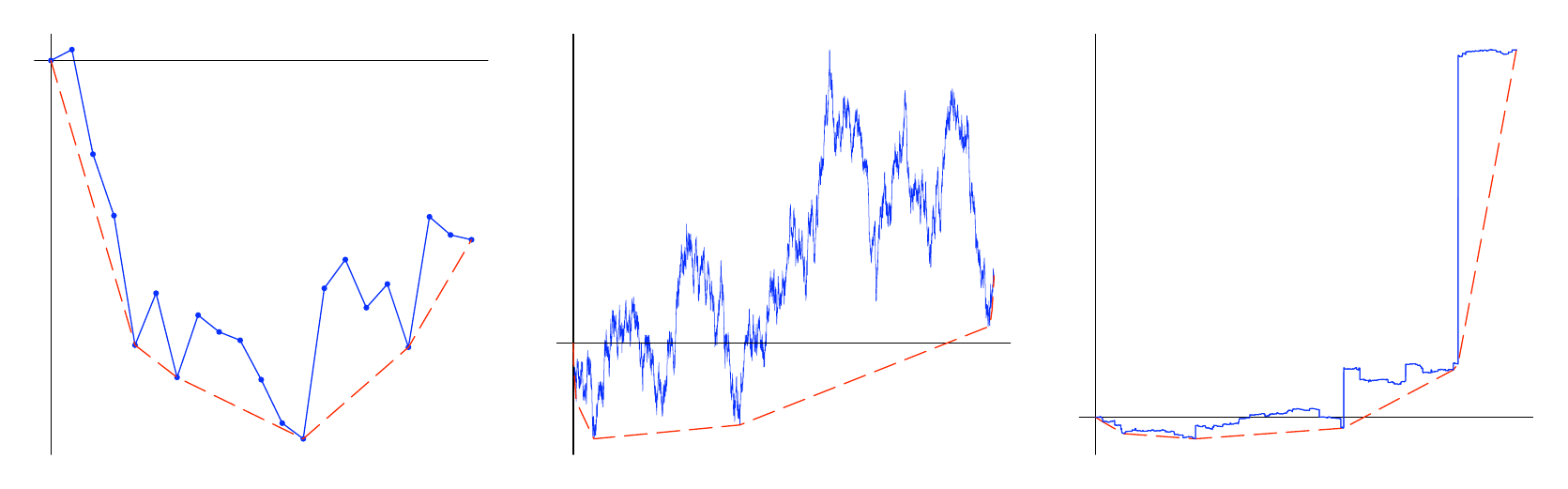}
\caption{Illustration of the convex minorant of a random walk, a Brownian motion, and a Cauchy process on a finite interval.}\label{1aa}
  \end{figure}

The {\em greatest convex minorant} (or simply \cmin\ for short) of a real-valued function $(x_u, u \in U)$ with domain $U$ contained in the real line
is the maximal convex function $(\unc_u, u \in I)$ defined on a closed interval $I$ containing $U$ with $\unc_u \le x_u$ for all $u \in U$.
A number of authors have provided descriptions of certain features of the \cmin\ for various stochastic processes such as random walks \cite{\cite{MR0058893b, MR0162302, MR994088,    MR0079851}},
Brownian motion \cite{MR2007793, brownistan, gboom83, p83, suidan01,MR770946}, Cauchy processes \cite{MR1747095}, Markov Processes \cite{Lachieze-Rey:2009fk}, and L\'evy processes (Chapter XI of \cite{MR1739699}).  Figure \ref{1aa} illustrates an instance of the \cmin\ 
for each of a random walk, a Brownian motion, and a Cauchy process on a finite interval.

The recent articles \cite{ap10,  piro10, pub10} provide a relatively complete description
of the \cmin\ of random walks, Brownian Motion, and L\'evy processes which not only encompass 
much of the somewhat ad hoc previous work on \cmin s, but also provide new tools to derive previously unknown properties of such \cmin s.  
These three articles together run well over 100 pages
and so the purpose of this note is to provide an overview of these works.  To this end, we will focus on stating results in
a streamlined fashion, referring to \cite{ap10,  piro10, pub10} where needed to furnish details and proofs.

The layout of the paper is as follows. In Section \ref{rw} we discuss the \cmin\ of a random walk following \cite{ap10}, and in Section \ref{lp} we
describe the limiting case of the results of Section \ref{rw} - the \cmin\ of a L\'evy process studied in \cite{pub10}.  In section \ref{bm} we provide
an overview of additional results in the special case of Brownian motion, found in \cite{piro10}. We conclude in section \ref{op} by presenting a selection of open problems.

\section{Random Walks}\label{rw}

Let $S_0 = 0$ and $S_j = \sum_{i=1}^j X_i$ for $1 \le j \le n$, where $X_1, \ldots , X_n$ 
are exchangeable random variables such that
almost surely no two subsets of $X_1,\ldots, X_n$ have the same arithmetic mean (satisfied for example if the  $X_i$
are i.i.d. with continuous distribution). Let $S^{[0,n]} := \{ (j,S_j) \: : \: 0 \le j \le n \}$, so that $S^{[0,n]}$ is the random walk of length $n$ with increments distributed
as $X_1,\ldots,X_n$.  As Figure \ref{1aa} indicates, the \cmin\ of $\son$ consists of piecewise linear segments which we refer to
as `faces.'  Let $F_n$ be the number of faces of $\unC^{[0,n]}$, the \cmin\ of $S^{[0,n]}$, and define
$$
0 < N_{n,1} < N_{n,1} + N_{n,2} < \cdots < N_{n,1} + \cdots + N_{n,F_n} = n
$$
be the successive indices $j$ with $0 \le j \le n$ such that $S_j=\unC_j$; we refer to $N_{n,i}$ as the `length'
of the $i$th face of $\unC^{[0,n]}$.  Finally, let $L_{n,1}, \ldots, L_{n,F_n}$ be the lengths of the faces
of $\unC^{[0,n]}$ arranged in non-decreasing order.  We refer to this sequence as the \emph{partition}
of $n$ generated by the \cmin\ of $\son$.  Recall the following classical result.
\begin{theorem}[\cite{MR0058893b, MR0162302, MR994088,    MR0079851}]
The sequence $L_{n,1}, \ldots, L_{n,F_n}$ of ranked lengths of faces of the \cmin\ of $\son$, a random walk with exchangeable
increments with almost surely no subset average ties has the same distribution as
the ranked cycle lengths of a uniformly chosen permutation of $n$ elements:
\[
\IP( F_n = k, L_i = n_i, 1 \le i \le k) = \prod_{j=1}^n \frac{1}{j^{a_j} a_j !}
\]
where $a_j := \# \{ i \: : \:  n_i = j \}$, and $n_1 \ge \ldots \ge n_k$ with $n_1 + \cdots + n_k = n$.
\end{theorem}

The following natural question was the starting point of our study of \cmin s.
\begin{quote}
Given the partition of $n$ generated by the faces of the \cmin\ of $\son$, how are the lengths ordered to form the
\emph{composition} of $n$ generated by the \cmin\ of $\son$?  
\end{quote}In the notation above, the sequence of variables $(N_{n,1}, \ldots, N_{n,F_n})$
is the composition of $n$ generated by the \cmin.

In the case that the $X_i$ are i.i.d. the answer to this question is easy to describe.  For $j=1,\ldots,n$ each face of length $j$
is assigned an increment distributed as $S_j$, independently of all other increments, and then the faces are ordered according
to increasing slope.  Formally, we have the following result.
\begin{theorem}[\label{ap1}\cite{ap10}]
Let $(N_{n,1}, \ldots, N_{n,F_n})$ be the composition of $n$ induced by the lengths  of the faces of the \cmin\ of $S^{[0,n]}$.
Assuming no subset average ties,
the joint distribution of $N_{n,1}, \ldots, N_{n,F_n}$ is given
by the formula
$$
\IP( F_n = k, N_{n,i} = n_i, 1 \le i \le k) = 
\IP \left( \frac { S_{n_1} ^{(1)} } { n_1 }  < \frac { S_{n_2} ^{(2)} } { n_2 }  <  \cdots < \frac { S_{n_k} ^{(k)}} { n_k } \right) \prod_{i = 1}^k \frac{1}{n_i}
$$
for all $n_1, \ldots, n_k$ with $n_1 + \cdots + n_k = n$, and
where for $1 \le i \le k$
$$
S_{n_i} ^{(i)}:=  S_{n_1 + \cdots + n_{i} } - S_{n_1 + \cdots + n_{i-1} } \ed S_{n_i}.
$$
In particular, if the $X_i$ are independent, then so are the $S_{n_j}^{(i)}$ for $1 \le i \le k$.
\end{theorem}

The special case of Cauchy increments gives rise to the following appealing version of Theorem \ref{ap1}.

\begin{corollary}
\label{cor:cauchy}
Suppose that the $X_i$ are independent and such that
$S_k/k$ has the same distribution for every $k$, as when
the $X_i$ have a Cauchy distribution. Then
$$
\IP( F_n = k; N_{n,i} = n_i, 1 \le i \le k) = 
\frac{1 }{k!} \prod_{i = 1}^k \frac{1}{n_i},
$$
and hence $\{ N_{n,i} \: : \: 1 \le i \le F_n \}$ has the same distribution as the composition of $n$ created by first choosing a random permutation of $n$ and then putting the cycle lengths in uniform random order.
\end{corollary}
Note that the continuum limit of this result can be read from Bertoin's work
\cite{MR1747095} and follows from the description provided in \cite{pub10} as discussed below.

\begin{figure}
  \centering
  \includegraphics[width=\textwidth]{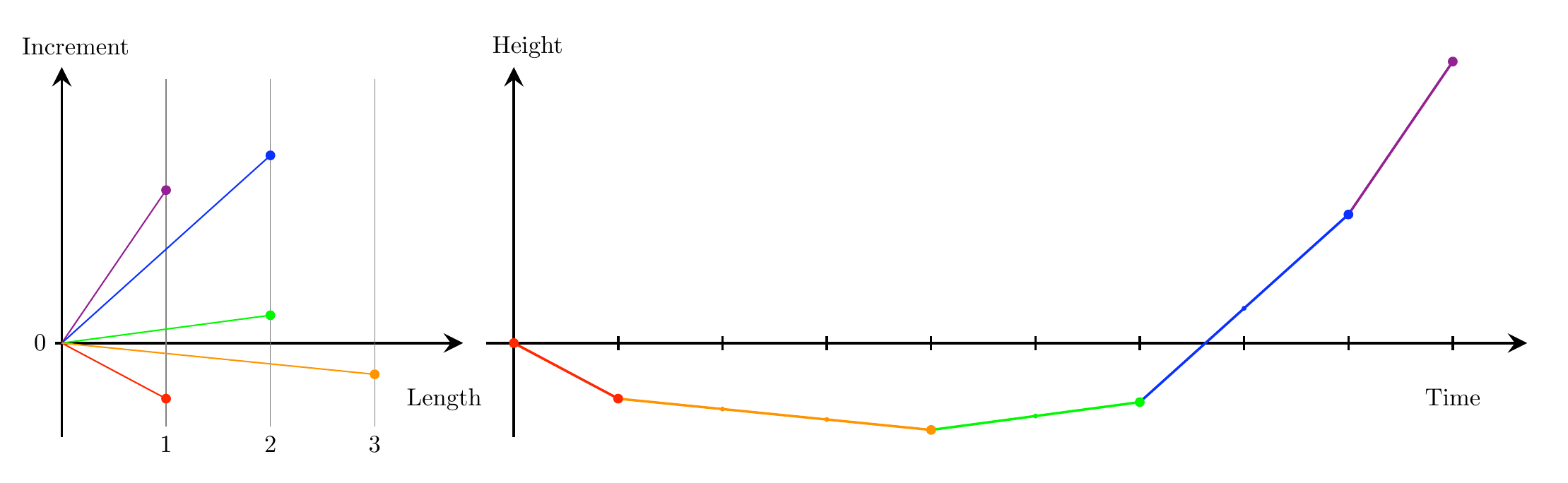}
\caption{Illustration of the convex minorant of a random walk as a point process.}\label{1b}
  \end{figure}

In order to proceed further, it is crucial that we introduce the representation of the \cmin\
as a point process of lengths and increments of the faces, where the lengths are chosen according
to the cycle structure of a random permutation of $n$ elements and the increments are chosen according to Theorem \ref{ap1} (independently if the $X_i$ are).  Figure
\ref{1b} illustrates this representation.

From this point, we can use Theorem \ref{ap1} to provide a construction of the \cmin\ of a random walk of a random length in the case of independent increments.
We already have some description in this case since we have a construction conditional on the length, but more can be said.
The work of Shepp and Lloyd \cite{MR0195117} on the
cycle structure of permutations combined with the forthcoming Proposition \ref{lem1} yield the following result.
\begin{theorem}[\cite{ap10}]\label{thm1}
Let $n(q)$ be a geometric random variable with parameter $1-q$; that is $\IP(n(q)\geq n)= q^n, n=0,1,\ldots$. If $X_1, X_2, \ldots $ are independent with common continuous distribution, then the point process of lengths and increments of faces the \cmin\ of $S^{[0,n(q)]}$ is a Poisson point process on $\{ 1,2,\ldots \} \times \mathbb{R}$ with intensity 
\ba
j^{-1} q^j \IP(S_j \in dx), \hspace{5mm} j = 1, 2, \ldots, \hspace{3mm} x \in \mathbb{R}. 
\ee
Moreover, let 
$T_i =  \sum_{l = 1}^i N_{n(q),l}$, $0 \le i \le F_{n(q)}$, be the consecutive indices at which $S^{[0,n(q)]}$ meets its \cmin, so that $T_0=0$ and $T_{F_{n(q)}} = n(q)$. Then 
the sequence of path segments
$$
\{ (S_{T_i + k} - S_{T_i},  0 \le k \le N_{n(q),i+1} ) , i = 0, \ldots , F_{n(q)} -1 \},
$$
is a list of the points of a Poisson point process in the space of finite random walk segments
$$
\{ (s_1, \ldots , s_j) \mbox{ for some } j = 1,2, \ldots \}
$$
whose intensity measure on paths of length $j$ is $q^j j^{-1}$ times the conditional distribution of the path $(S_1, \ldots , S_j)$ given that $S_k > (k/j) S_j $ for all $1 \le k \le j-1 $.
\end{theorem}

An important facet of the Poisson point process description is that it provides a decomposition of
a random walk up to the index of its minimum.  For example, the description of Theorem \ref{thm1}
is a more complete description of the \cmin\ of a random walk which was the basis for
Spitzer's combinatorial identity \cite{MR0079851}. 
\begin{theorem}[\cite{MR0079851}]\label{spit}
Let $X_1, X_2, \ldots $ be independent with common continuous distribution, $S_0=0$,
 $S_k=\sum_{i=1}^{k}X_i$ for $k\geq1$, and $M_n := \min_{0 \le k \le n} S_k$.  Then
 \ba
\sum_{n = 0}^\infty q^{n} \E e^{i t M_n } = \exp \left( \sum_{k = 1}^\infty \frac{q^k}{k} \E e^{it S_k^-}\right),
\ee
where $S_k^-=\min\{S_k, 0\}$.
\end{theorem}

%
Now, by letting $q$ tend to one in Theorem \ref{thm1}, we obtain a description of the \cmin\ of $S^{[0,\infty)}$, a random
walk on $[0,\infty)$.
\begin{theorem}[\cite{ap10}]
If $X_1, X_2, \ldots $ are independent with common continuous distribution with $\IE X_1 \in (-\infty,\infty]$, then the point process of lengths and increments of faces the \cmin\ of $S^{[0,\infty)}$ is a Poisson point process on $\{ 1,2,\ldots \} \times \mathbb{R}$ with intensity 
\ba
j^{-1} \IP(S_j \in dx), \hspace{5mm} j = 1, 2, \ldots, \hspace{3mm} \frac{x}{j}<\IE X_1. 
\ee
\end{theorem}
Similar to Theorem \ref{thm1}, there is a companion path space statement which we omit for the sake of brevity.

The key to the results above is a certain property of a transformation of the walk $\son$, which
not only yields the results above, but also provides a construction of the walk jointly with its \cmin.
We will call this transformation the `3214' transformation, as it is described by first dividing the walk $\son$
into four consecutive paths, and then reordering these four pieces with the third one first, the second one second, and so on.

\begin{figure}
  \centering
  \includegraphics[width=.49\textwidth]{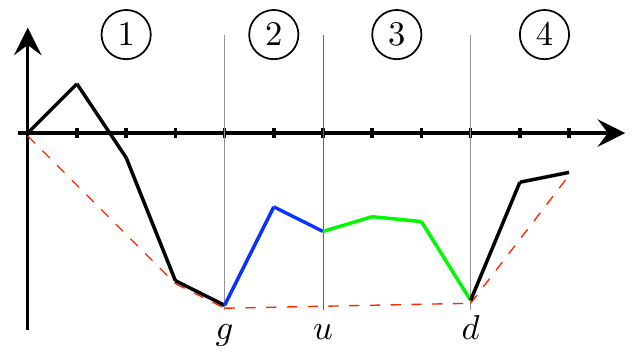}
   \includegraphics[width=.49\textwidth]{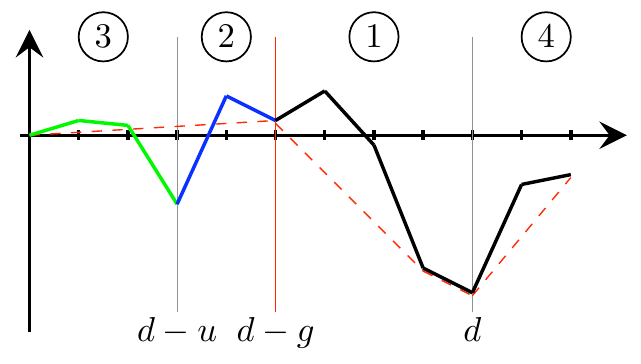}
\caption{Notation and application of the `3214' transformation.}\label{2a}
  \end{figure}

The `3214' transform of $\son$ is generated by a random variable $U$ which is uniform on
$\{1, \ldots, n\}$ and is independent of $\son$.  Given $U=u$, we then define $g$ and $d$ as the indices of the 
left and right endpoints of the face of the \cmin\ of $\son$ straddling the index $u$.  Note that $g$ and $d$ are almost surely
well defined by this description due to the no subset average ties assumption.  Consider the four paths of the random
walk on the intervals $[0,g]$, $[g,u]$, $[u,d]$, and $[d,n]$.   
With this setup, the `3214' transform is defined by reordering the four path fragments of $\son$
described above to form a new walk path $\son_U$ in the order $3-2-1-4$. 
Figure \ref{2a} illustrates the
notation and provides an example of the transformation.

The following lemma summarizes the crucial feature of this transform.

\begin{lemma}[\cite{ap10}]\label{3214}
Let $X_1, \ldots, X_n$ be exchangeable random variables with no subset average ties and $\son$ the random
walk generated by the $X_i$.  Let $U$ uniform on $\{1,\ldots,n\}$, independent of the $X_i$ and $\son_U$ the `3214' transform of $\son$
generated by $U$.  If $g$ and $d$ are the indices of the endpoints of the face of the \cmin\ of $\son$
to the left and right of $U$, then
\ba
\left(U,\son\right)\ed \left(d-g, \son_U\right).
\ee
\end{lemma}

To see how Lemma \ref{3214} corroborates the story above, we introduce discrete \emph{uniform stick breaking}
on $[0,n]$, one of the many well-known representations of the distribution of the cycle lengths of a uniformly chosen
permutation on $n$ elements.

\begin{definition}
For an integer $n$, define the discrete uniform stick breaking sequence of random variables $M_{n,1}, \ldots, M_{n,K_n}$ as follows.
\begin{itemize}
\item $M_{n,1}$ is uniform on $\{1, \ldots, n\}$.
\item For $i\geq1$, if $\sum_{j=1}^iM_{n,j}<n$, then $M_{n,i+1}$ is uniform on $\left\{1, n- \sum_{j=1}^iM_{n,j}\right\}$.
\item For $i\geq1$, if $\sum_{j=1}^iM_{n,j}=n$, then set $K_n=i$, and end the process.
\end{itemize}
We refer to the variables  $L_{n,1},\ldots, L_{n,K_n}$ defined to be $M_{n,1}, \ldots, M_{n,K_n}$ rearranged in non-increasing order
as the \emph{partition} of $n$ generated by uniform stick breaking.
\end{definition}
To be explicit, we state the following well-known proposition (see \cite{MR2245368}).
\begin{proposition}
\label{lem1}
The partition of $n$ generated by uniform stick breaking has the same distribution as
the ranked cycle lengths of a uniformly chosen permutation of $n$ elements.
\end{proposition}

From this point, some consideration yields the following implications of Lemma \ref{3214} for a walk with i.i.d. increments and no subset average ties:
\begin{itemize}
\item The lengths of the faces of the \cmin\ of $\son$ are distributed as discrete stick breaking.
\item Conditional on the lengths of the faces of the \cmin\ of $\son$, the excursions above the segments are independent.
\item Given a segment of length $j$, the excursion above the segment can be realized as the unique cyclic permutation of a random
walk of length $j$ equal in distribution to  $S^{[0,j]}$ which yields a \cmin\ of exactly one segment.  
\end{itemize}
The last item is similar in spirit to Vervaat's transform of a Brownian bridge to an excursion \cite{MR515820}.  As 
this transformation is not well developed for random walks and L\'evy processes in general (some statements for L\'evy processes are
found in \cite{gub11}), this last item
carries real content.

\begin{figure}
  \centering
  \includegraphics[scale=1.3]{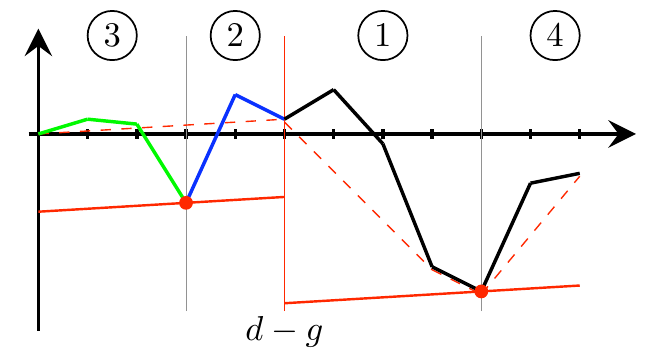}
\caption{Illustration of the inverse `3214' transform for the walk of Figure \ref{2a}.}\label{2c}
  \end{figure}

The proof 
of Lemma \ref{3214} essentially follows from two observations.  The first
is that given the values of the increments $X_1=x_1$, \ldots, $X_n=x_n$, $S_j$
is distributed as $\sum_{i=1}^jx_{\sigma_i}$ for $j=1, \ldots, n$ and where $\sigma$ is a permutation chosen uniformly at random.
From this point we only need to show that the `3214' transformation is a bijection between $\{1,\ldots,n\}\times$ `paths generated from
permutations of $x_1, \ldots, x_n$' for fixed increments $x_i$ having no subset average ties.  The bijection is easily verified after 
noting that for a given value of $d-g$, the indices at which Segment 1 meets Segment 4 and Segment 3 meets Segment 2 are found by
raising a line with slope equal to the mean of the first $d-g$ increments.  Figure \ref{2c} illustrates this inverse transformation.

If we remove the assumption that almost surely, no two subsets of the $X_i$ have the same mean, then the process of generating excursions described above may generate excursions that meet the corresponding face of the convex minorant at points other than the end points, and excursions that have the same slope. This implies that there is
not necessarily a unique cyclic permutation transforming a walk into an excursion, and neither is there necessarily a unique ordering of the excursions that puts them in non-decreasing order of slope.  Such technical issues can be dealt with in
a straightforward manner, but the details, found in  \cite{ap10}, are a little laborious.

\section{L\'evy processes}\label{lp}

A real valued process $X$ is a L\'evy process on $[0,\infty)$ if
$X_0 = 0$, $X$ is cadlag (right continuous with left limits), and  $X$ has independent and stationary increments.
As is well known, L\'evy processes are the continuous scaling limits of discrete time random walks generated by i.i.d. increments, so
it is not surprising that continuous analogs of the results of Section \ref{rw} hold for L\'evy processes.
However, there are a few interesting wrinkles not present in the discrete case and many technical
details to be considered in pushing the discrete results to the limit.  We restrict our analysis
to the case that $X_t$ has continuous distribution for all $t>0$, which is equivalent to 
the assumption that $X$ is not a compound Poisson process with drift.

In analogy to the case of random walk, we can view the intervals that a L\'evy process is strictly greater 
than its \cmin\ on $[0,1]$ as an interval partition of the unit interval.  
The formal statement of this last fact is proved in \cite{pub10}, but also
intersects with the work \cite{Lachieze-Rey:2009fk}.

\begin{proposition}[\cite{pub10}]\label{prop1}
Let $X$ be a L\'evy process with continuous distribution  and $\unC$  the convex minorant of $X$ on $[0,1]$.  Let
$\mathcal{O}=\{s\in(0,t):\unC_s<X_s\wedge X_{s-}\}$ and  $\mathcal{I}$ be the set of connected components of $\mathcal{O}$.
The following conditions hold almost surely:
\begin{enumerate}
\item\label{ComplementIsLightProperty} The open set $\mathcal{O}=\{s\in(0,t):\unC_s<X_s\wedge X_{s-}\}$ has Lebesgue measure $1$.
\item $\mathcal{I}$ is a set of disjoint intervals and the closure of its union is $[0,1]$.
\item\label{SimplicityOfSlopesProperty}If $(g_1,d_1)$ and $(g_2,d_2)$ are distinct intervals of $\mathcal{I}$, then the
slopes of $\unC$ over those intervals  differ:
\ba
\frac{\unC_{d_1}-\unC_{g_1}}{d_1-g_1}\neq\frac{\unC_{d_2}-\unC_{g_2}}{d_2-g_2}.
\ee
\end{enumerate}
\end{proposition}

For each $(g,d)\in\mathcal{I}$, we refer to $g$ and
$d$ as \emph{vertices}, the \emph{length} is $d-g$, the \emph{increment} is $\unC_d-\unC_g$, and the \emph{slope} is
$(\unC_d-\unC_g)/(d-g)$.

%
%

Because the partition of $n$ generated by the \cmin\ of an i.i.d. generated random walk of $n$ steps is distributed
as the partition of $n$ generated by the cycles of a random permutation for any increment distribution, we might hope that
a similar universal result holds for L\'evy processes and also that this universal result might be a limiting
continuous distribution of the cycle structure of a random permutation.  This is indeed the case, but
before stating our result we define this continuous limit.

\begin{definition}
Define the continuous uniform stick breaking sequence of random variables as the sequence  $L_{1}, L_2, \ldots$ defined as follows.
\begin{itemize}
\item $L_{1}$ is uniform on $[0,1]$.
\item For $i\geq1$, $L_{i+1}$ is uniform on $\left[0, 1-\sum_{j=1}^iL_{j}\right]$.
\end{itemize}
We refer to the variables  $L_{1}, L_2, \ldots, $ rearranged in non-increasing order
as the \emph{partition} of $[0,1]$ generated by uniform stick breaking.
\end{definition}

The variables $L_1, L_2, \ldots$ almost surely sum to one and their law once arranged in decreasing order is referred to as the 
\emph{Poisson-Dirichlet distribution with parameter one} which is the limiting distribution
of the cycle structure of a permutation chosen uniformly at random (see \cite{MR2245368}).
We can now state the following result and note that
a proof in the special case of Brownian motion was sketched in \cite{suidan01}.

\begin{theorem}[\cite{pub10}]\label{thm4}
The sequence of ranked lengths of faces of the \cmin\ of a L\'evy process with continuous distributions
has the Poisson-Dirichlet distribution with parameter one.
\end{theorem}

In light of Theorem \ref{thm4}, the following natural question arises.  Given the interval partition of $[0,1]$ generated by
the \cmin\ of a L\'evy process $X$ with continuous distribution, how are the intervals ordered to form the interval \emph{composition}
of $[0,1]$ generated by the \cmin\ of $X$?

In total analogy with the answer for i.i.d. random walks, the answer to this question is easy to describe.  Given the interval of
length $l$, the increment is distributed as $X_l$ independent of all other increments, and the faces are ordered according to
increasing slope. 
\begin{theorem}[\cite{pub10}]\label{thmppl}
Let $X$ be a L\'evy process with continuous distribution and let $\{(g_i,d_i)\}_{i\geq1}$ denote the intervals 
of $\unC$, the \cmin\ of $X$ on $[0,1]$.  Let $L_1, L_2, \ldots$ be generated by uniform stick breaking,
$S_0:=0$, and for $i\geq1$, define $S_i:=\sum_{j=1}^iL_i$.  Then we have the following equality in
distribution between sequences:
\ba
\left(\left(d_i-g_i,\unC_{d_i}-\unC_{g_i}\right), i\geq 1\right) \ed \left(\left(L_i, X_{S_i}-X_{S_{i-1}}\right),i\geq 1\right).
\ee
\end{theorem}
We note here that applying the theorem to a Cauchy process yields the main result of Bertoin \cite{MR1747095}
and also shows the composition 
generated by the \cmin\ on $[0,1]$ is a uniform ordering of the generated partition; c.f. Corollary \ref{cor:cauchy}.

We can also consider the \cmin\ of a L\'evy process $X$ on an interval of a random exponential length independent of $X$  to obtain
the following analog of Theorem \ref{thm1}.
\begin{theorem}[\cite{pub10}]\label{thmpl}
Let $T$ be a rate $\theta$ exponential random variable, $X$ a L\'evy
process with continuous distribution which is independent of $T$ and let $\unC^T$ denote the \cmin\ of $X$ on $[0,T]$.
The point process 
\ba
\left\{\left(d_i-g_i,\unC^T_{d_i}-\unC^T_{g_i}\right), i\geq 1\right\}
\ee
generated by the lengths $t$ and increments $x$ of $\unC^T$ 
has the same distribution as a Poisson point process with intensity measure
\ba
\mu(dt,dx)=\frac{e^{-\theta t}}{t}\, dt\IP(X_t\in dx).
\ee
\end{theorem}

By integrating out the independent exponential variable, we can also use Theorem \ref{thmpl}
gain insight into the structure of the \cmin\ of a L\'evy process on $[0,1]$.  

For example, this program yields the following neat dichotomy for stable L\'evy processes.
\begin{proposition}[\cite{pub10}]\label{thmst}
Let $X$ be a symmetric stable process with parameter $\alpha$; that is $\IE e^{iuX_t}=e^{-t|u|^\alpha}$, and let
$\mathcal{S}$ be the set of slopes and $\mathcal{T}$ be the set of times of vertices of the \cmin\ of $X$ on $[0,1]$.  
\begin{itemize}
\item If $1< \alpha \leq 2$, then $\mathcal{S}$ has no accumulation points, $\mathcal{S}\cap(a,\infty)$ and $\mathcal{S}\cap(-\infty, -a)$ 
are infinite for $a>0$, and $\mathcal{T}$ has accumulation points at zero and one only.
\item If $0<\alpha \leq 1$, then $\mathcal{S}$ is dense $\IR$ and every point of $\mathcal{T}$ is an accumulation point. 
\end{itemize}
\end{proposition}

By letting $\theta$ tend to zero in Theorem \ref{thmpl},
we obtain a description of the \cmin\ of $X$ on $[0,\infty)$ which was also derived in \cite{MR1739699}.
\begin{theorem}[\cite{MR1739699, pub10}]
If $X$ is a L\'evy process with continuous distribution and and 
\ba
I:=\lim\inf_{t \rightarrow \infty}\frac{X_t}{t}\in(-\infty,\infty],
\ee
then
the lengths $t$ and increments $x$ of the \cmin\ of $X$ on $[0,\infty)$ is a Poisson point process with intensity
\ba
\frac{\IP(X_t\in dx)}{t}, \hspace{5mm} x<I t.
\ee
\end{theorem}
Both of the previous theorems carry an It\^o type excursion theory analogous to that of Theorem \ref{thm1} for random walks, see \cite[Thm. 4]{pub10}.

Theorems \ref{thm4} and \ref{thmppl}
follow from a direct analog of Lemma \ref{3214} for a `3214' transform for L\'evy processes.
The proof uses limiting arguments which crucially hinge on certain regularity conditions for L\'evy processes governing
the behavior of the process at the vertices of the \cmin.


\section{Brownian Motion}\label{bm}

Since Brownian motion is a L\'evy process (and stable with index $2$), the results of the previous section
apply to the \cmin\ of Brownian motion, and some of these results were known (from \cite{MR2007793, brownistan, gboom83, p83, suidan01}).
However, Brownian motion offers extra analysis due to its special properties among L\'evy processes (e.g. continuity and time
inversion).  We begin by noting the following special case of Theorem \ref{thmpl}.
\begin{theorem}[\cite{gboom83}]\label{bmpp}
Let $\Gamma_1$ be an exponential random variable with rate one.
The lengths~$x$ and slopes~$s$ of the faces of the \cmin\ of a Brownian motion on 
$[0, \Gamma_1]$
form a Poisson point process on $\IR^+ \times \IR$ with intensity measure
\ba
\frac{\exp\{-\frac{x}{2}\left(2+s^{2}\right)\}}  {\sqrt{2 \pi x}} \, ds\,dx, \hspace{5mm} x\geq0, s\in \IR. 
\ee
\end{theorem}

As with random walks and L\'evy processes, the minimum on $[0,T]$ of a Brownian motion is a distinguished point of the \cmin\
and the process after the minimum can be described by restricting the point process of slopes and increments to
those points with positive slopes.  Due to Proposition \ref{thmst}, we can define
\ba
\alpha_0 < \alpha_1 < \alpha_2 < \cdots < 1 
\ee
with $\alpha_n \uparrow 1$ as $n \te \infty$ to be times
of vertices of the convex minorant of a Brownian motion $B$ on $[0,1]$, arranged relative to
\ba
\alpha_0 := \argmin_{ 0 \leq t \leq 1 } B_t. 
\ee
Brownian scaling and Theorem \ref{bmpp} yield an implicit description of the distribution of the sequence $(\alpha_i)_{i\geq0}$.
Moreover, Denisov's decomposition for Brownian motion at the minimum \cite{MR726906} implies that
the process after the minimum is a Brownian meander, for which we now provide an alternate description.
First we make the following definition.
\begin{definition}
We say that a sequence of random variables $(\tau_n, \rho_n)_{n\geq0}$ satisfies the \tr\ recursion if
for all $n \ge 0$:
\ba
\rho_{n+1} = U_n \rho_{n}  \mbox{ \, and \, } \tau_{n+1} = \frac{ \tau_n  \rho_{n+1}^2 }{ \tau_n Z_{n+1}^2 + \rho_{n+1}^2 }
\ee
for the two independent sequences of i.i.d. uniform $(0,1)$ variables $U_n$ and i.i.d. squares of standard normal random variables $Z_n^2$, both independent of $(\tau_0, \rho_0)$.
\end{definition}

\begin{figure}
\begin{center}
\includegraphics[width=.5\textwidth]{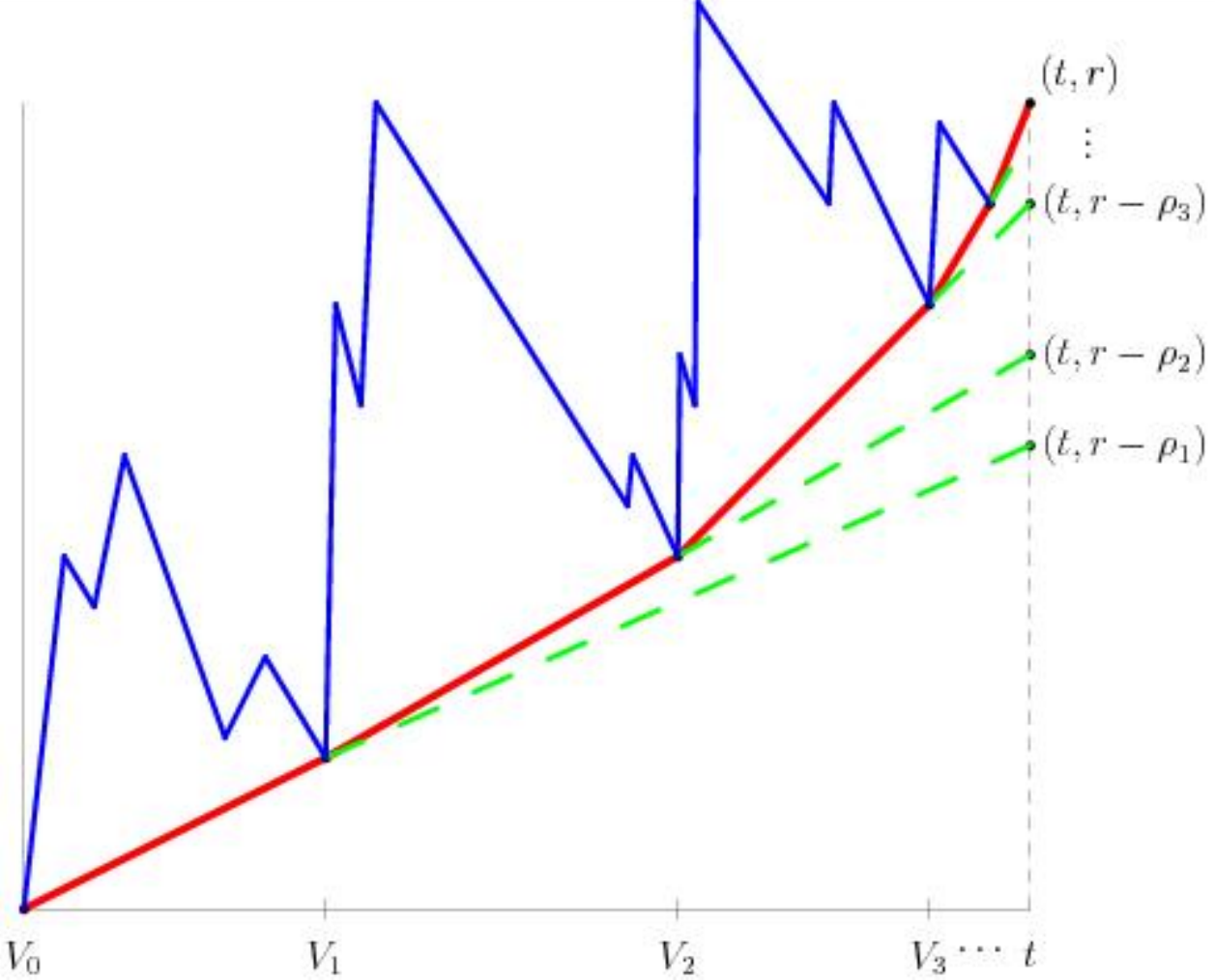}
\end{center}
\caption{An illustration of the notation of Theorem \ref{mean}.  The blue line represents a Brownian meander of length $t$,
and the red line its \cmin.  Note also that $V_i:= t-\tau_i$ for $i=0, 1, \ldots$}\label{22}
\end{figure}

\begin{theorem}[\cite{piro10}]\label{mean}
Let $(X(v), 0 \le v \le t)$ be a Brownian meander of length $t$, and let
$(\unC(v), 0 \le v \le t)$ be its \cmin. The vertices of $(\unC(v), 0 \le v \le t)$ occur at times
$0 = V_0 < V_1 < V_2  < \cdots$ with $\lim_n V_n = t$. Let $\tau_n:= t - V_n$ so $\tau_0 = t > \tau_1 > \tau_2 > \cdots$
with $\lim_n  \tau_n  = 0$. Let $\rho_0 = X(t)$ and for $n\geq1$ let $\rho_0- \rho_n$ denote the intercept at time $t$ of the line extending the segment of the
\cmin\ of $X$ on the interval $(V_{n-1}, V_n)$. The \cmin\ $\unC$ of $X$ is uniquely determined by the sequence
of pairs $(\tau_n, \rho_n)$ for $n = 1, 2, \ldots$ which satisfies the \tr\ recursion with
\ba
\rho_0 \ed \sqrt{2 t \Gamma_1 } \mbox{ and } \tau_0 = t,
\ee
where $\Gamma_1$ is an exponential random variable with rate one.
\end{theorem}

In \cite{piro10}, we use the descriptions provided by Theorems \ref{bmpp} and \ref{mean} (and the interplay between them
as per Denisov's decomposition \cite{MR726906})
to derive various properties about the \cmin\ of Brownian motion on $[0,1]$, such as formulas for
densities of the $\alpha_i$.  We also use the equivalence of the two descriptions to
discover new identities between related quantities in each description.
We conclude this section with an elementary example of such an identity; we leave a direct proof as a challenge to the reader.
\begin{corollary}[\cite{piro10}]\label{bewid}
Let $W$ and $Z$ be standard normal random variables, $U$ uniform on $(0,1)$, and $R$ Rayleigh distributed having density $re^{-r^2/2}$, $r>0$.  If
all of these variables are independent, then 
\ba
\left(\frac{ W^2 + (1-U) ^2 R^2 }{ 1 + U^2 R^2/Z^2 }, \frac{(1-U) R}{\sqrt{T}}\right) \ed \left(Z^2, \frac{(1-U) R}{\sqrt{T}}\right).
\ee
Note that the two coordinate variables on the right are independent.
\end{corollary}

\section{Open Problems}
\label{op}
We end this note with a list of open problems.
\begin{itemize}
\item Under what conditions is the right derivative of the \cmin\ of a L\'evy process with continuous distribution discrete, continuous, or mixed?
\item Provide a description of the \cmin\ of a continuous time process with \emph{exchangeable} increments. 
\item Provide a framework independent of the \cmin\ of Brownian motion that explains the equivalence of the Poisson point process
of Theorem \ref{bmpp} with the sequential description of Theorem \ref{mean}. 
\item Is there a version of the sequential description of Theorem \ref{mean} for random walks or L\'evy processes?
\end{itemize}
\bibliographystyle{amsplain}
\bibliography{convex_minorant}

\providecommand{\bysame}{\leavevmode\hbox to3em{\hrulefill}\thinspace}
\providecommand{\MR}{\relax\ifhmode\unskip\space\fi MR }
\providecommand{\MRhref}[2]{%
  \href{http://www.ams.org/mathscinet-getitem?mr=#1}{#2}
}
\providecommand{\href}[2]{#2}
\begin{thebibliography}{10}

\bibitem{ap10}
J.~Abramson and J.~Pitman, \emph{Concave majorants of random walks and related
  {P}oisson processes}, \url{http://arxiv.org/abs/1011.3262}, 2010.

\bibitem{MR0058893b}
Erik~Sparre Andersen, \emph{On the fluctuations of sums of random variables
  {II}}, Math. Scand. \textbf{2} (1954), 195--223.

\bibitem{MR770946}
R.~F. Bass, \emph{Markov processes and convex minorants}, Seminar on
  probability, {XVIII}, Lecture Notes in Math., vol. 1059, Springer, Berlin,
  1984, pp.~29--41. \MR{MR770946 (86d:60086)}

\bibitem{MR1747095}
J.~Bertoin, \emph{The convex minorant of the {C}auchy process}, Electron. Comm.
  Probab. \textbf{5} (2000), 51--55 (electronic). \MR{MR1747095 (2001i:60081)}

\bibitem{MR0162302}
H.~D. Brunk, \emph{A generalization of {S}pitzer's combinatorial lemma}, Z.
  Wahrscheinlichkeitstheorie und Verw. Gebiete \textbf{2} (1964), 395--405
  (1964). \MR{MR0162302 (28 \#5501)}

\bibitem{MR2007793}
C.~Carolan and R.~Dykstra, \emph{Characterization of the least concave majorant
  of {B}rownian motion, conditional on a vertex point, with application to
  construction}, Ann. Inst. Statist. Math. \textbf{55} (2003), no.~3, 487--497.
  \MR{MR2007793 (2004h:60123)}

\bibitem{brownistan}
E.~{\c{C}}inlar, \emph{Sunset over {B}rownistan}, Stochastic Process. Appl.
  \textbf{40} (1992), no.~1, 45--53. \MR{MR1145458 (93f:60120)}

\bibitem{MR726906}
I.~V. Denisov, \emph{Random walk and the {W}iener process considered from a
  maximum point}, Teor. Veroyatnost. i Primenen. \textbf{28} (1983), no.~4,
  785--788. \MR{MR726906 (85f:60117)}

\bibitem{MR994088}
C.~M. Goldie, \emph{Records, permutations and greatest convex minorants}, Math.
  Proc. Cambridge Philos. Soc. \textbf{106} (1989), no.~1, 169--177.
  \MR{MR994088 (91a:60138)}

\bibitem{gboom83}
P.~Groeneboom, \emph{The concave majorant of {B}rownian motion}, Ann. Probab.
  \textbf{11} (1983), no.~4, 1016--1027. \MR{MR714964 (85h:60119)}

\bibitem{Lachieze-Rey:2009fk}
Raphael Lachieze-Rey, \emph{Concave majorant of stochastic processes and
  {B}urgers turbulence}, \url{http://arxiv.org/abs/0909.1088v2}, 2009.

\bibitem{MR1739699}
M.~Nagasawa, \emph{Stochastic processes in quantum physics}, Monographs in
  Mathematics, vol.~94, Birkh\"auser Verlag, Basel, 2000. \MR{MR1739699
  (2001g:60148)}

\bibitem{p83}
J.~Pitman, \emph{Remarks on the convex minorant of {B}rownian motion}, Seminar
  on Stochastic Processes, 1982, Birkh{\"a}user, Boston, 1983, pp.~219--227.
  \MR{MR733673}

\bibitem{MR2245368}
\bysame, \emph{Combinatorial stochastic processes}, Lecture Notes in
  Mathematics, vol. 1875, Springer-Verlag, Berlin, 2006, Lectures from the 32nd
  Summer School on Probability Theory held in Saint-Flour, July 7--24, 2002,
  With a foreword by Jean Picard. \MR{MR2245368 (2008c:60001)}

\bibitem{piro10}
J.~Pitman and N.~Ross, \emph{The convex minorant of {B}rownian motion, meander,
  and bridge.}, \url{http://arxiv.org/abs/1011.3073}, 2010.

\bibitem{pub10}
J.~Pitman and G.~Uribe~Bravo, \emph{The convex minorant of a {L}{\'e}vy
  process}, \url{http://arxiv.org/abs/1011.3069}, 2010.

\bibitem{MR0195117}
L.~A. Shepp and S.~P. Lloyd, \emph{Ordered cycle lengths in a random
  permutation}, Trans. Amer. Math. Soc. \textbf{121} (1966), 340--357.
  \MR{0195117 (33 \#3320)}

\bibitem{MR0079851}
F.~Spitzer, \emph{A combinatorial lemma and its application to probability
  theory}, Trans. Amer. Math. Soc. \textbf{82} (1956), 323--339. \MR{MR0079851
  (18,156e)}

\bibitem{suidan01}
T.~M. Suidan, \emph{Convex minorants of random walks and {B}rownian motion},
  Teor. Veroyatnost. i Primenen. \textbf{46} (2001), no.~3, 498--512.
  \MR{MR1978665 (2004d:60097)}

\bibitem{gub11}
G.~Uribe~Bravo, \emph{Bridges of {L}\'evy processes conditioned to stay
  positive.}, \url{http://arxiv.org/abs/1101.4184}, 2011.

\bibitem{MR515820}
Wim Vervaat, \emph{A relation between {B}rownian bridge and {B}rownian
  excursion}, Ann. Probab. \textbf{7} (1979), no.~1, 143--149. \MR{515820
  (80b:60107)}

\end{thebibliography}
\end{document}